\documentclass[12pt]{article}
\usepackage{amsthm,amsfonts,amsmath, amscd, bbold, bm}
\usepackage{mathrsfs}
\usepackage{amssymb} 
\usepackage{dsfont}
\usepackage{accents}

\usepackage{mathrsfs}
\usepackage{amscd}
\usepackage[active]{srcltx}
\usepackage{verbatim}
\usepackage{xcolor}

\newtheorem{theorem}{Theorem}[section]
\newtheorem{lemma}[theorem]{Lemma}

\newtheorem{remark}[theorem]{Remark}
\newtheorem{definition}[theorem]{Definition}

\newcommand{\C}{{\mathord{\mathbb C}}}
\newcommand{\Z}{{\mathord{\mathbb Z}}}
\newcommand{\N}{{\mathord{\mathbb N}}}
\newcommand{\E}{{\mathord{\mathbb E}}}

\newcommand{\tr}{{\mathop{\rm Tr}}}

\newcommand{\cQ}{{\mathord{\mathscr Q}}}
\newcommand{\cL}{{\mathord{\mathscr L}}}
\newcommand{\cK}{{\mathord{\mathscr K}}}

\newcommand{\cP}{{\mathord{\mathscr P}}}

\newcommand{\V}{{\mathcal{V}}}
\newcommand{\G}{{\mathcal{G}}}
\renewcommand{\H}{{\mathcal H}}
\renewcommand{\V}{{\mathcal V}}
\newcommand{\B}{{\mathcal B}}

\renewcommand{\G}{{\mathcal G}}

\renewcommand{\E}{{\mathcal E}}

\newcommand{\bma}{{\bm{\alpha}}}
\newcommand{\bmb}{{\bm{\beta}}}
\newcommand{\bmg}{{\bm{\gamma}}}
\newcommand{\bmd}{{\bm{\delta}}}

\newcommand{\bmet}{{\bm{\eta}}}
\newcommand{\bmz}{{\bm{\zeta}}}

\newcommand{\one}{{\mathds{1}}}

\newcommand{\Dens}{{\mathfrak S}}

\begin{document}

\title{Spectrum for some  Quantum Markov semigroups describing $N$-particle systems evolving under a binary collision mechanism}

\author{\vspace{5pt} Eric A. Carlen$^1$   and
Michael P. Loss$^{2}$ \\
\vspace{5pt}\small{$1.$ Department of Mathematics, Hill Center,}\\[-6pt]
\small{Rutgers University,
110 Frelinghuysen Road
Piscataway NJ 08854-8019 USA}\\
\vspace{5pt}\small{$2.$ School of Mathematics,
Georgia Tech } \\[-6pt]
\small{Atlanta GA 80332
  }\\
 }
\date{\today}
\maketitle 


\let\thefootnote\relax\footnote{
\copyright \, 2023 by the authors. This paper may be
reproduced, in its entirety, for non-commercial purposes. Partial support by US National Science Foundation DMS-2055282 (E.A.C.) and  DMS-2154340 (M.L.) is acknowledged.}

\begin{abstract}
We compute the spectrum for a class of quantum Markov semigroups describing systems of $N$  particle interacting through a binary collision mechanism.
These quantum Markov semgroups are associated to a novel kind of quantum random walk on graphs, with the graph structure arising naturally 
in the quantization of the classical Kac model, and we show that the spectrum of the generator of the quantum Markov semigroup is closely 
related to the spectrum of the Laplacian on the corresponding graph.   For the direct analog of the original classical Kac model, 
we determine the exact spectral gap for 
the quantum generator.   We also give a new and simple method for studying the spectrum of certain graph Laplacians. 
\end{abstract}

\medskip
\leftline{\footnotesize{\qquad Mathematics subject
classification numbers: 05C50, 60J10}}
\leftline{\footnotesize{\qquad Key Words: Quantum Markov semigroup, graph Laplacian}}

\section{Introduction} \label{intro}

We study the rate of approach to equilibrium in a quantum version of the classical Kac model that was developed in \cite{CCL19}. 
The original Kac model \cite{K56,K59}  concerns a a dilute gas of $N$ molecules interacting through pair-wise collision that conserve the energy, and in this model, the collision mechanism gives rise to a jump process on a continuous state space, the ``energy sphere'' of the $N$ particles, and  the Kolmogorov forward equation for this jump process is known as the {\em Kac Master Equation}.

Interactions between molecules are properly described by quantum mechanics, and the model introduced in  \cite{CCL19} is a natural adaptation of 
Kac's classical model to the quantum setting, in which the Kac Master Equation becomes an equation of Lindblad type, the 
{\em Quantum Kac Master Equation} (QKME).
The  assumptions on the collision mechanism in \cite{CCL19} were rather general, and the paper  concentrated on general features such as classifying the
equilibrium states, proving propagation of chaos (see \cite{CCL19}) and studying the resulting non-linear quantum Boltzmann equation. Little was said about the actual evolution
and about the spectral properties of the generator of the quantum Markov semigroup, which is the subject of this paper.

In the present work we consider  a special collision rule that is described below and which is the direct analog of the one first considered by Kac.   This allows a much more detailed analysis. 
As we show, the quantization leads to a graph structure on an orthonormal basis of eigenstates of the $N$-particle energy operator. We shall give
a complete description of the spectrum of the QKME generator
 in terms of spectra of the graph Laplacian on this graph.  This yields a description of the quantum evolution in terms of 
 the eigenvectors of these graph Laplacians. 

In fact, one may view the QKME as a quantum Kolmogorov forward equation for a sort of quantum random walk on the 
graph associated to the collision mechanism.  However, this sort of quantum random walk is different from the the class of random 
walks on graphs introduced by 
Aharonov, Ambainis, Kempe and Vazirani \cite{AAKV}.  In particular, the graph structures emerges naturally from the dynamics, and is not present from the beginning. Of course, it is not surprising to see discrete mathematical structures emerging from quantization. 

In the model we focus on here, the graphs that arise  turn out to be the closely connected with {\em multislice} (see \cite{FOW19}), a 
natural generalization of the sliced Boolean cube.  Recently proved ergodic properties of the random walk on slices of the multislice  
\cite{Fil16,FOW19,S20} then become relevant to QKME.  However, the Kac model perspective sheds light on spectral analysis of Laplacians on graphs: We give a simple proof of the result of Caputo, Liggett and Richthammer \cite{CLR} that the spectral gap of the Laplacian on the multislice for $N$ particles has the same value, namely $N$, for all non-trivial connected components of the graph. This leads to an exact determination of the spectral gap for our quantum Kac model.

\subsection{Description of the quantum Kac model}

Let us recall the setting in \cite{CCL19}. Consider a $d$-dimensional Hilbert space $\mathcal{H}$ and a single particle Hamiltonian $h$ acting on $\H$  with $d$  eigenvalues $e_1, \dots, e_d$.  For the moment, we impose no conditions on the spectrum. Let 
$\{\psi_1, \dots, \psi_d\}$ be an orthonormal basis for $\H$ with $h\psi_j = e_j\psi_j$ for all $1 \leq j \leq d$.  

Define the $N$-particle Hamiltonian $H_N$  on $\H_N = \otimes^N\H$  as a sum of $N$ terms
$$
H_N = (h\otimes \one \otimes \cdots \otimes \one) + \cdots  + (\one \otimes \cdots \otimes h \otimes \cdots \otimes \one) + \cdots + (\one\otimes \cdots \otimes \one \otimes h)\ ,
$$
where $\one$ is the identity on $\H$, and where in the $n$th term, there is  a single $h$ which  is in the $n$th  position. 

 
 The eigenvalues of  $H_N$ are indexed by the multi-indices $\bma = (\alpha_1,\dots,\alpha_N) \in \{1,\dots,d\}^N$, and are  given by 
\begin{equation}\label{ealpha}
e(\bma) = e_{\alpha_1} + \cdots + e_{\alpha_N} \ .
\end{equation}
Defining
$$
\Psi_\bma := \psi_{\alpha_1} \otimes \cdots \otimes \psi_{\alpha_N} \ ,
$$
$\{ \Psi_\bma \ :\  \bma \in \V_N\}$ is an orthonormal basis of $\H_N$ consisting of eigenvectors of $H_N$.
For a multi-index $\bma$, and $k \in \V$ define the `occupation numbers'
$$
{\bf k}(\bma) = (k_1(\bma), \dots,k_d(\bma))
\quad,\qquad
k_j(\bma) = |\{ 1\leq m \leq N \:\  \alpha_m = j\ \}|
$$
where for a set $A$,  $|A|$ denotes the cardinality of $A$.  Thus, an alternate to the formula \eqref{ealpha} for $e(\bma)$ is 
\begin{equation}\label{ealpha2}
e(\bma) = \sum_{j=1}^{d} k_j(\bma)e_j\ .
\end{equation}
We also define, for all $1 \leq m < n \le N$,
\begin{equation}\label{emna}
e_{m,n}(\bma) = e_{\alpha_m} + e_{\alpha_m}\ .
\end{equation}

 We will be investigating a quantum kinetic time evolution in which particles interact through instantaneous
``pair collisions'', and between collisions, particles do not interact, and the energy of the system  is then given by $H_N$, which includes no interactions.   

Consider a pair collision between particles $m$ and $n$ with $1 \leq m < n \leq N$.  Suppose initially that the state of the system is $\Psi_\bma$.
Then after the collision, the state of the system will be a linear combination of state of the form $\Psi_\bmb$ where
\begin{equation}\label{adjac1}
\beta_p = \alpha_p \quad{\rm for}\quad p \neq m,n
\end{equation}
since only particles $m$ and $n$ are involved in the collision, and such that
\begin{equation}\label{adjac2}
e_{m,n}(\bma) = e_{m,n}(\bmb)\ ,
\end{equation}
since the collisions conserve energy.

\subsection{The Quantum Kac Master Equation} 

We now describe a  quantum analog of the Kac Master Equation \cite{K56} that will govern the time evolution of our system. 
The evolution of the system is driven by a Poisson stream of binary collisions that conserve energy.
To describe  binary collisions  in our $N$ particle system, we first consider the case $N=2$, and will then ``lift'' our 
constructions to the full $N$-particle model.  Let $U$ be some unitary operator on $\H_2$, and define 
define $U_{1,2}$
to be the unitary operator on $\H_N$ given by
$$U_{1,2} (\phi_1\otimes \phi_2\otimes\phi_3\otimes  \cdots \otimes \phi_N) = U(\phi_1 \otimes \phi_2)\otimes \phi_3 \otimes \cdots \otimes \phi_N\ .$$
In the same way, for $1\leq m< n \leq N$, we define $U_{m,n}$ so that it acts on the $m$th and $n$th factors of $\H$. The unitary $U_{m,n}$ 
describes the effects of a particular type of collision between particles $m$ and $n$.   Since the collisions conserve energy, we require that each 
$U_{m,n}$ commutes with $H_N$, and clearly this is the same thing as requiring that $U$ commutes with $H_2$.   
We regard all unitaries $U$ that commute with $H_2$ as describing a kinematically possible collision.

We discuss here a simple collision model, one that is an analog of the original classical Kac model, in which we average uniformly over all 
unitaries $U$ that commute with $H_2$.   In this simplest quantum Kac model, the quantum operation that describes a 
collision between two particles, viewed as acting on ${\mathcal B}({\mathcal H}_2)$,  is given by
\begin{equation}\label{Qspec}
\cQ(X) = \int_{\mathcal U} U^* XU {\rm d}\nu\ ,
\end{equation}
where ${\mathcal U}$ is the group of all unitaries commuting with $H_2$ and $\nu$ is normalized Haar measure on it.  
Evidently, $\cQ$ is completely positive, unital and trace preserving. Since  Haar measure on ${\mathcal U}$ is invariant under the map the map $U \mapsto U^*$, 
$\cQ$ is self adjoint with respect to the Hilbert-Schmidt inner product on $\B(\H_2)$.

It is easy to see that $\cQ$ is a projection, and in fact $\cQ(X)$ is the ``pinching'' operation that projects onto the $C^*$ subalgebra of  
${\mathcal B}({\mathcal H}_2)$ consisting of functions of $H_2$. Therefore, if  
$$
H_2 = \sum_{E \in \Sigma_2} E P_E
$$
is the spectral decomposition of $H_2$,
\begin{equation}\label{Qdef2}
\cQ(X)  =  \sum_{E \in \Sigma_2} \frac{1}{D(E)} \tr[ XP_E]P_E\ ,
\end{equation}
where $D(E)= \tr[P_E]$ is the degeneracy of the eigenvalue $E$.
From here, it is easy to write down a Kraus representation of $\cQ$ in terms of the  orthonormal basis $\{\psi_{\alpha_1}\otimes \psi_{\alpha_2}\ :\ 1 \leq \alpha_1,\alpha_2 \leq d\}$ of $\H_2$.  Since
$$\tr[XP_E] = \sum_{e_{\alpha_1}+e_{\alpha_2}= E}\langle \psi_{\alpha_1}\otimes \psi_{\alpha_2} ,X\psi_{\alpha_1}\otimes \psi_{\alpha_2}\rangle \quad{\rm and}\quad P_E := 
\sum_{e_{\beta_1}+e_{\beta_2}= E}|  \psi_{\beta_1}\otimes \psi_{\beta_2} \rangle\langle\psi_{\beta_1}\otimes \psi_{\beta_2} |\ ,
$$if we define
$$
F_{\alpha_1\alpha_2;\beta_1\beta_2} := |\psi_{\alpha_1}\otimes \psi_{\alpha_2} \rangle\langle\psi_{\beta_1}\otimes \psi_{\beta_2} |\ ,
$$
then
\begin{equation}\label{Qdef2}
\cQ(X)  =  \sum_{E \in \Sigma_2} \frac{1}{D(E)} \left( \sum_{e_{\alpha_1}+e_{\alpha_2} =  e_{\beta_1}+e_{\beta_2}= E} 
F_{\alpha_1\alpha_2;\beta_1\beta_2}^*X F_{\alpha_1\alpha_2;\beta_1\beta_2}\right)\ .
\end{equation}


We now lift  $\cQ$ to $\H_N$ as described in the beginning of this section by letting fixing $1 \leq m < n\leq N$ and replacing each $U$ by $U_{m,n}$, 
and thus obtaining $\cQ_{m,n}$, which now describes the averaged effect of a collision between particles $m$ and $n$.    For every $1 \leq m < n \leq N$,  
$\cQ_{m,n}$ is completely positive,  self-adjoint with respect to the Hilbert-Schmidt inner product on $\B(\H_N)$, unital and trace 
preserving because it inherits these properties form $\cQ$.  

In lifting $\cQ$ up to $\H_N$ as described, the trace in \eqref{Qdef2} becomes the partial trace over the $m$th and $n$th  factors in $\H_N$.  
For $(m,n)=(1,2)$ one easily finds
$$
\cQ_{1,2}(X) =\sum_{\alpha_1, \alpha_2,\beta_1,\beta_2 ; e_{\alpha_1}+e_{\alpha_2} = e_{\beta_1}+e_{\beta_2}} \frac{1}{D(e_{\alpha_1} + e_{\alpha_2})}  E_{\alpha_1\alpha_2;\beta_1\beta_2}^*X E_{\alpha_1\alpha_2;\beta_1 \beta_2}
$$
where
$$
E_{\alpha_1\alpha_2;\beta_1 \beta_2} =  |\psi_{\alpha_1}\otimes\psi_{\alpha_2}\rangle \langle \psi_{\beta_1}\otimes\psi_{\beta_2}| \otimes \one_{N-2} =
 |\psi_{\alpha_1}\rangle \langle \psi_{\beta_1}| \otimes |\psi_{\alpha_2}\rangle \langle \psi_{\beta_2}| \otimes \one_{N-2}  \ ,
$$
and $\one_{N-2}$ is the identity operator on the remaining factor. For general $m<n$ the rank one operators $|\psi_{\alpha_1}\rangle \langle \psi_{\beta_1}|$ and 
$ |\psi_{\alpha_2}\rangle \langle \psi_{\beta_2}|$
 should be inserted as factors $m$ respectively $n$  in the tensor product and the remaining factors are the identity matrices. Therefore, for $1 \leq m < n \leq N$, 
 \begin{equation}\label{kraus}
\cQ_{m,n}(X) =\sum_{\alpha_m, \alpha_n,\beta_m,\beta_n ; e_{\alpha_m}+e_{\alpha_n} = e_{\beta_m}+e_{\beta_n}} \frac{1}{D(e_{\alpha_m} + e_{\alpha_n})}  E_{\alpha_m\alpha_n;\beta_m\beta_n}^*X E_{\alpha_m\alpha_n;\beta_m \beta_n}
\end{equation}
where
$$
E_{\alpha_m\alpha_n;\beta_m \beta_n} =   |\psi_{\alpha_m}\rangle \langle \psi_{\beta_m}| \otimes  |\psi_{\alpha_n}\rangle \langle \psi_{\beta_n}|\otimes \one_{N-2} \ ,
$$
with the subscripts indicating the factors on which the first two terms operate.  Finally define
$$
\cQ(X)=  {\binom{N}{2}}^{-1}\sum_{m<n} \cQ_{mn}(X) \ .
$$
Note that \eqref{kraus} gives the Kraus form of the completely positive operator $\cQ_{m,n}$, and then summing, we have the Kraus 
form of $\cQ$.  From this one can easily write down the Lindblad form of the generator $\cL_N$, but here it turns out to be more 
convenient to work directly with \eqref{kraus}.

We recall the following definitions from \cite{CCL19}:

\begin{definition}
Define the operators $\cQ_N$ and $\cL_N$ on $\mathcal{B}(\H_N)$ by
\begin{equation}\label{quankac9}
\cQ_N  =  {\binom{N}{2}}^{-1}\sum_{m<n} \cQ_{m,n} \quad{\rm and}\quad 
\cL_N  = N(\cQ_N - {\mathds{1}}_{\H_N})\ .
\end{equation} 
\end{definition}

$\cQ_N$ is completely positive,  self-adjoint with respect to the Hilbert-Schmidt inner product on $\B(\H_N)$, unital and trace 
preserving because the set of such operators is convex.  Note that \eqref{kraus} gives the Kraus form \cite{K71} of the completely positive operator $\cQ_{m,n}$, and then summing, we have the Kraus 
form of $\cQ$.  From this one can easily write down the Lindblad form \cite{L76} of the generator $\cL_N$, but here it turns out to be more 
convenient to work directly with \eqref{kraus}.

\begin{definition}
The {\em Quantum Kac Master Equation} (QKME) is the evolution equation on $\Dens(\H_N)$ given by
\begin{equation}\label{quankac12}
\frac{{\rm d}}{{\rm d} t} \varrho(t)  = \cL_N \varrho(t) \ .
\end{equation}

Since $\|\cK_N\|_\infty \leq 2N$, the QKME is solved by exponentiation:
For each $t\geq 0$ ,we may define an operator $\cP_{N,t}$ on each $\B(\H_N)$ by 
\begin{equation}\label{quankac11}
\cP_{N,t} A  = \sum_{k=1}^\infty e^{-Nt}\frac{(Nt)^k}{k!} \cQ^k_N A  = e^{t\cL_N} A\ .
\end{equation} 
Then the unique solution $\varrho(t)$ of the QKME satisfying $\varrho(0) = \varrho_0\in \Dens(\H_N)$ is
$\varrho(t) = \cP_{N,t}\varrho_0$. 
\end{definition}

\section{The graph structure induced by the collision rules}

The pair collision rules \eqref{adjac1} and \eqref{adjac2} induce a graph structure on $\{1,\dots,d\}^N$ which then becomes its set of vertices.  This graph structure is fundamental to our analysis of the spectral properties of the QKME generator.

\begin{definition}Let $\V_N$ denote $\{1,\dots,d\}^N$, the set of multi-indices, now considered as vertices in a graph. Two vertices $\bma, \bmb \in \V_N$ are {\em  adjacent}   in case $\bma \neq \bmb$ and \eqref{adjac1} and \eqref{adjac2} are satisfied.
Any such pair of adjacent vertices defines an edge, and we denote the set of all such edges by ${\mathcal E}_N$.  If $\bmg$ is adjacent to $\bmd$, we write $[\bmg,\bmd]$ to denote the corresponding edge. Note that $[\bmg,\bmd] = [\bmd,\bmg]$.  We denote this graph by $\mathcal{G}_N$.  
\end{definition}

It is evident that if $[\bmg,\bmd]  \in {\mathcal E}_N$, then $e(\bmg) = e(\bmd)$, but the converse need not be true,  as we explain below.

\begin{remark} If the eigenvalues $\{e_1,\dots,e_d\}$ are all distinct, and if $\bma$ and $\bmb$ are adjacent, they differ in exactly two places, so that to each edge $[\alpha,\beta]\in \E_N$, there is a uniquely determined pair $(m,n)$, $1 \leq m < n \leq N$ such that $\alpha_m\neq \beta_m$ and $\alpha_n\neq \beta_n$.   If the eigenvalues of $h$ are degenerate, then there will exist 
adjacent $\bma$ and $\bmb$ that differ in a single index so that for some $m$, $\alpha_m \neq \beta_m$, but $e_{\alpha_m} = e_{\beta_m}$, and hence \eqref{adjac1} and \eqref{adjac2} are satisfied with $\alpha_j = \beta_j$ for all $j\neq m$, and any choice of $n\neq m$. 
Such a transition can be viewed as resulting from a collision between particle $m$ and any other particle $n$ that has 
somehow ``catalyzed'' the transition of particle $m$.  Thus, 
no single pair collision can be identified with the edge $[\bma,\bmb]$. {\em Going forward, we assume the eigenvalues of $h$ are non-degenerate.}
\end{remark}

\begin{definition} Let  $\sigma(H_2)$ denote the spectrum of $H_2$.  For each $E\in \sigma(H_2)$, recall that  $D(E)$ is the dimension of the corresponding eigenspace.  As remarked above, under the assumption that the eigenvalues of $h$ are non-degenerate, to each $[\bma,\bmb]\in \E_N$, there is a unique
$(m,n)$, $1 \leq m < n \leq N$ such that $\alpha_m\neq \beta_m$ and $\alpha_n\neq \beta_n$.  Since $\bma$ and $\bmb$ are adjacent, \eqref{adjac2} is satisfied. Define
$$
D[\bma,\bmb] = D(e_{\alpha_m} + e_{\alpha_n})  = D(e_{\beta_m} + e_{\beta_n})\ 
$$
for this uniquely determined pair $(m,n)$. 
 \end{definition}

\begin{definition} For each $\bma\in \V_N$, let $v(\bma)$ denote the valency of $\bma$; i.e., the number of vertices in $\G_N$ that are adjacent to $\bma$. 
\end{definition}

For each $\bma$, there are $D(e_{\alpha_m} +e_{\alpha_n})-1$ vertices that are adjacent to $\bma$ through collision involving particles $m$ and $n$.  Summing over $m< n$, the valency of $\bma$ is given by
$$
v(\bma) = \sum_{m< n} (D(e_{\alpha_m} +e_{\alpha_n})-1)\ ,
$$
or alternatively,
\begin{equation}\label{valency0}
v(\bma) = \sum_{i< j}^d k_i(\bma)k_j(\bma)(D(e_i+e_j) -1)
 \end{equation}
Under the non-degeneracy assumption for $i< j$, $D(e_i+e_j) \geq 2$, and hence
\begin{equation}\label{valency}
 v(\bma) \geq \sum_{i < j} k_i(\bma)k_j(\bma)\ .
 \end{equation}
 
 Given a finite undirected graph $\G$ with vertex set $\V$ and edge set $\E$,  the {\em graph Laplacian} $\Delta_\G$ is the operator on functions $f$ on $\V$ given by
\begin{equation}\label{Lapdef}
\Delta_{\G}f(x) = \sum_{y \in \V \ :\ \{x,y\} \in \E} (f(x) - f(y))\ .
\end{equation}
Note that
\begin{equation}\label{Lapdef2}
\sum_{x\in \V} f(x) \Delta_{\G}f(x) =\sum_{y \in \V \ :\ \{x,y\} \in \E} \frac12(f(x) - f(y))^2\ .
\end{equation}
This computation shows that if $\mu_{\V}$ denotes the uniform probability measure on $\V$, 
$L_\G$ is a positive semi-definite operator on $L^2(\mu_{\V})$, and that the constant function $f(x) = 1$ for all $x\in \V$ is an 
eigenfunction with eigenvalue $0$.  (This is the standard sign convention for the Laplacian in graph theory.) 

It follows from \eqref{Lapdef2} that $ \Delta_{\G}f =0$ if and only if $f$ is constant on each 
connected component of $\G$. Hence, on a connected graph $\G$, $0$ is an eigenvalue of multiplicity one, and the eigenspace 
is spanned by the constant vector. The quadratic form on the right in \eqref{Lapdef2} is called the {\em Dirichlet form} of the graph Laplacian.

\begin{definition} Let $\G$ be a finite connected graph with vertex set $\V$. The {\em spectral gap} of $\G$, $\Gamma_\G$,  
is the least non-zero eigenvalue of $\Delta_\G$.
\end{definition}

By the Rayleigh-Ritz variational principle,
\begin{equation}\label{RRV}
\Gamma_\G = \inf\left \{ \int_{\V} f(x) \Delta_{\G} f(x){\rm d}{\mu_{\V}} \ :\ \int_{\V} f(x) {\rm d}{\mu_{\V}} = 0\ ,\  \int_{\V} |f(x)|^2 {\rm d}{\mu_{\V}} =1\ \right\}\ .
\end{equation}

\section{The spectrum of the Kac generator $\cL_N$}

When the spectrum of the single particle hamiltonian $h$ satisfies a strong non-degeneracy condition, it is possible to give a simple and complete spectral decomposition of the generator $\cL_N$ for the QKME.

\begin{definition} The spectrum of $h$ is {\em strongly non-degenerate} in case:

\smallskip
\noindent{\it(1)} 
The spectrum of $h$  is such that the spectrum of $H_2$  is non-degenerate on the symmetric subspace of $\H_2 =\H\otimes \H$.  That is,  for any $1 \leq j_1,j_2,j_3,j_4 \leq d$,
 \begin{equation*}
 e_{j_1} + e_{j_2} = e_{j_3} + e_{j_4}  \iff \{j_1,j_2\} = \{j_3,j_4\}\ .
 \end{equation*}

\smallskip
\noindent{\it(2)} For each $N \geq 3$ and each $E\in \sigma(H_N)$, 
the pair of equations
\begin{eqnarray*}
 \sum_{m=1}^{d} k_m e_n =  E\nonumber \quad{\rm and}\quad  \sum_{m=1}^{d} k_m  &=&  N
  \end{eqnarray*}
 has exactly one solution for each $E$ in the spectrum of $\H_N$.   
\end{definition}

Note that condition {\it (2)} is satisfied if the $\{e_1,\dots,e_d\}$ are linearly independent over the rational numbers, and in this sense the condition is generically satisfied.

We shall show that when $h$ is strongly non-degenerate, $\B(\H_N)$ has a decomposition into a direct sum of subspaces, each of which is invariant under the Kac generator $\cL_N$, and moreover, the restriction of $-\cL_N$ to each of these invariant subspaces is unitarily equivalent to $\Delta_{\G_M}$ plus an explicit non-negative  constant for some $M\leq N$.

\subsection{Direct sum decomposition of $\B(\H_N)$}

First,  we fix some notation. For $\bma,\bmb\in \V_N$, define 
$$
F_{\bma\bmb} := |\Psi_\bma\rangle\langle \Psi_\bmb|\ .
$$
Since $\{\Psi_\bma\}_{\bma\in \V_N}$ is an orthonormal basis for $\H_N$,  $\{\ F_{\bma\bmb}\ :\ \bma,\bmb\in \V_N\ \}$ is an orthonormal basis of $\B(\H_N)$. 
 
 \begin{definition}\label{CoincDef}  For any  $\bmd,\bmg\in \V_N$, define $C_{\bmd\bmg}\subseteq\{1,\dots,N\}$ to be the {\em coincidence set} of this pair of vertices. That is,
$$
C_{\bmd\bmg} = \{ \ j\ :\ \delta_j = \gamma_j\ \}\ .
$$
 
 Next, let $S\subsetneq \{1,\dots,N\}$:
 
 \smallskip
 \noindent{\it (1)}  Define the set of {\em exterior pair configurations} $E_S$ to be the set of pairs
$\bmz,\bmet\in \{1,\dots,d\}^{S^c}$ such that for all $j\in S^c$, $\zeta_j\neq \eta_j$.

\smallskip
 \noindent{\it (2)} For each $S\subsetneq \{1,\dots,N\}$ and each $(\bmz,\bmet)\in E_S$, define
 $\B_{S;\bmz,\bmet}$ to be the subspace of $\B(\H_N)$
 spanned by the operators $F_{\bmd,\bmg}$ such that $C_{\bmd\bmg}= S$, and such that for $j\in S^c$, $\delta_j = \zeta_j$ and $ \gamma_j= \eta_j$\ .
 
 \smallskip
 \noindent{\it (3)} For $S=\{1,\dots,N\}$, define $\mathcal{C}_N = {\rm span}(\{ F_{\bmg,\bmg}\ :\ \bmg\in \V_N\ \})$. That is, $\mathcal{C}_N$ is the span of the $F_{\bmd,\bmd}$ such that $C_{\bmd\bmg} = \{1,\dots,\N\}$ in which case there are no exterior configurations to be considered.  Note that $\mathcal{C}_N$ is not only a subspace of $\B(\H_N)$, it is a commutative subalgebra of $\B(\H_N)$. It is called the {\em classical subalgebra of $\B(\H_N)$.}
 \end{definition}
 
 \begin{lemma}\label{DecompLem} The $N$-particle space $H_N$ is the direct sum of the subspaces $\B_{S;\bmz,\bmet}$; that is
 \begin{equation}\label{DirSumInv}
 \H_N = \mathcal{C}_N \oplus \left(\bigoplus_{S\subsetneq\{1,\dots,N\}\ ,\ (\bmz,\bmet)\in E_S} \B_{S;\bmz,\bmet}\right)\ .
 \end{equation}
 \end{lemma}
 
 \begin{proof} It is evident that for each $\bmd,\bmg\in \V_N$, $F_{\bmd,\bmg}$ belongs to $\B_{S;\bmz,\bmet}$ if and only if $S= C_{\bmd\bmg}$ and for each $j \notin S$,$\delta_j = \zeta_j$ and $\gamma_j = \eta_j$.  Thus, each of the basis vectors $F_{\bmd,\bmg}$ to exactly one of the spaces $\B_{J;\bmz,\bmet}$, and this proves
 \eqref{DirSumInv}. 
%
 \end{proof} 
 
 For $M\in \N$, let $L^2(\G_M)$ denote the Hilbert space obtained by equipping $\G_M$ with counting measure. For each multi-index $\bma = (\alpha_1,\dots,\alpha_M)$,
define a function $g_\bmb$ on $\V_M$ by 
$g_\bmb(\bma) = \delta_{\bma,\bmb}$. Then $\{g_\bmb\}_{\bma\in \V_M}$ is an orthonormal basis for $L^2(\G_M)$.

 Fix  some $S\subsetneq\{1,\dots,N\}$ with cardinality $|S| = M$,   and write it in the form $\{j_1,\dots,j_M\}$. Define the map $k_S:S\to \{1,\dots,M\}$ by $k_S(j_k) = k$.
In addition, for $M < N$, fix some
$(\bmz,\bmet)\in E_S$. Then
define a unitary map $U_{S;\bmz\bmet}$ from  $L^2(\G_M)$  to $\B_{S;\bmz,\bmet}$ by linearly extending
\begin{equation}\label{UnitDef}
U_{S;\bmz\bmet}(g_{\bmb})= F_{\bmg,\bmd}\quad{\rm where}\quad  \begin{cases} \gamma_j = \delta_j = \beta_{k_S(j)}& j\in S\\
\gamma_j = \zeta_j\ ,\ \delta_j = \eta_j & j\notin S\end{cases}\ .
\end{equation}
This map is unitary because it takes an orthonormal basis of $L^2(\G_M)$  to an orthonormal basis of $\B_{S;\bmz,\bmet}$. For $S= \{1,\dots,N\}$ there is the simpler unitary map $U_{\mathcal{C}_N}$ from $L^2(\G_M)$ onto $\mathcal{C}_N$ defined by  $U_{\mathcal{C}_N}(g_\bmb)= F_{\bmb\bmb}$.

\begin{theorem}\label{GrLem2} Let $S\subsetneq\{1,\dots,N\}$ with $|S| = M$, and let $(\bmz,\bmet)\in E_S$. 
Define $r := N-M$.  For all $f\in L^2(\G_{M})$, and all $\bma,\in \V_{M}$,
\begin{equation}\label{GraphCon}
U_{S;\bmz\bmet}^*\left(-\cL_N \right) U_{S;\bmz\bmet}f(\bma) 
 = \frac{1}{N-1} \Delta_{\G_M} f(\bma) + \left(2r\frac{N}{N-1} \textcolor{red} {- }\frac{r(r+1)}{N-1}\right)f(\bma)\ .
\end{equation}
for $S= \{1,\dots,N\}$ we have 
\begin{equation}\label{GraphCon22}
U_{\mathcal{C}_N}^*\left(-\cL_N \right) U_{\mathcal{C}_N} 
 = \frac{2}{N-1}\Delta_{\G_N} f(\bma) \ .
\end{equation}
\end{theorem}

Theorem~\ref{GrLem2} reduces the study of the spectrum of $\cL_N$ to the study of the spectrum of the graph Laplacians $\Delta_{\G_M}$, $M \leq N$.  
The rest of this section is devoted to the proof of Theorem~\ref{GrLem2}.  In the next section, we turn to the analysis of the spectrum of 
$\Delta_{\G_M}$, $M \leq N$. We give a new a relatively simple proof of some known relevant facts about the spectrum, and in particular, 
the determination of the spectral gap.

\begin{lemma}  Suppose that $h$ is strongly non-degenerate. If $\bmg$ and $\bmd$ are adjacent in $\G_N$, then $\bmg$ and $\bmd$
differ differ by a pair transposition. The valency of $\bma\in \G_N$ is given by
\begin{equation}\label{valency2}
 v(\bma) = \sum_{i < j} k_i(\bma)k_j(\bma)\ .
 \end{equation}

\end{lemma} 

\begin{proof} Since $\bmg$ and $\bmd$ are adjacent, they differ in exactly 2 indices, say $m$ and $n$. But then since $k_j(\bmd) = k_j(\bmg)$ for each $j=1,\dots,d$, 
$\{\gamma_m,\gamma_n\} = \{\delta_m,\delta_n\}$. Then since  $\bmg\neq\bmd$, $\delta_m = \gamma_n$  and $\delta_n = \gamma_m$. Thus $\bmg$ and $\bmd$
differ  by a pair transposition. Finally, \eqref{valency2} follows from \eqref{valency0} since for $i < j$, $D(e_i+e_j) =2$. 
\end{proof}

 The operator
$\cQ$ turns out to have a fairly simple matrix representation in the $\{F_{\bma,\bmb}\}$ basis, as we now show.   It will be useful to define
$$
f_{\alpha\beta} = |\psi_\alpha\rangle\langle \psi_\beta|\ .
$$
We shall also make use of the {\em swap map}, both as a map from $\V_N$ into itself, and as a unitary operator on $\H_N$. First,  a map form $\V_N$ into itself, for
$m < n$, and $\bma\in \V_N$, $S_{mn}(\bma)$ is defined by
$$
(S_{mn}(\bma))_j = \begin{cases}\bma_n & j =m\\ \bma_m &j=n\\ \bma_j& j\neq m,n\end{cases}\ .
$$
The {\it swap operator} is the unitary operator on $\H_N$ defined by
$$
S_{mn} \Psi_\alpha = \Psi_{S_{mn}(\bma)}\ .$$
\begin{lemma}\label{WalkLem}
We have that
\begin{eqnarray}
{\binom{N}{2}}\cQ(F_{\bmg\bmd} ) &=&\frac 12\sum_{m<n} \delta_{\gamma_m\delta_m} \delta_{\gamma_n\delta_n} 
 [F_{\bmg\bmd} + S_{mn}F_{\bmg\bmd}S_{mn}]\nonumber \\
&=&\frac 12\sum_{m<n} \delta_{\gamma_m\delta_m} \delta_{\gamma_n\delta_n} 
 [F_{\bmg\bmd} + F_{S_{mn}(\bmg)S_{mn}(\bmd)}]  \ .\label{WalkForm}
\end{eqnarray}
\end{lemma}
\begin{proof}
 For simplicity we pick the pair $(m,n)= (1,2)$. Then
$$
 E_{\alpha_1\alpha_2;\beta_1\beta_2}^*F_{\bmg\bmd}  E_{\alpha_1\alpha_2;\beta_1 \beta_2}
$$
$$
=[f_{\beta_1\alpha_1}\otimes f_{\beta_2\alpha_2}\otimes I_{N-2}] [\otimes_{j=1}^N f_{\gamma_j \delta_j}][ f_{\alpha_1\beta_1}\otimes f_{\alpha_2\beta_2}\otimes I_{N-2}]
$$
$$
=f_{\beta_1\alpha_1}f_{\gamma_1 \delta_1}  f_{\alpha_1\beta_1} \otimes f_{\beta_2\alpha_2} f_{\gamma_2 \delta_2}f_{\alpha_2\beta_2} \otimes_{j=3}^Nf_{\gamma_j \delta_j}
$$
$$
=\delta_{\alpha_1\gamma_1} \delta_{\delta_1\alpha_1} \delta_{\alpha_2\gamma_2} \delta_{\delta_2\alpha_2}  f_{\beta_1\beta_1}\otimes f_{\beta_2\beta_2}  \otimes_{j=3}^Nf_{\gamma_j \delta_j} \ .
$$
The sum
$$
\sum_{\alpha_1, \alpha_2,\beta_1,\beta_2 ; e_{\alpha_1}+e_{\alpha_2} = e_{\beta_1}+e_{\beta_2}} \frac{1}{D(e_{\alpha_1} + e_{\alpha_2})}  E_{\alpha_1\alpha_2;\beta_1\beta_2}^*F_{\bmg\bmd} E_{\alpha_1\alpha_2;\beta_1 \beta_2}
$$
has a contribution from the non-degenerate eigenvalues of $H_2$ given by $2e, e\in \sigma(h)$
$$
 E_{\alpha_1\alpha_1;\alpha_1\alpha_1}^*F_{\bmg\bmd} E_{\alpha_1\alpha_1;\alpha_1 \alpha_1}
=\delta_{\gamma_1\delta_1}\delta_{\gamma_2\delta_2} \delta_{\gamma_1\gamma_2}  f_{\gamma_1\gamma_1}\otimes f_{\gamma_1\gamma_1}  \otimes_{j=3}^Nf_{\gamma_j \delta_j} = 
\delta_{\gamma_1\delta_1}\delta_{\gamma_2\delta_2} \delta_{\gamma_1\gamma_2}F_{\bmg\bmd} \ .
$$
and a contribution from the doubly degenerate eigenvalues of $h$, $e_i+e_j, i\not= j$ 
$$
\frac 12 \sum_{\alpha_1, \alpha_2,\beta_1,\beta_2 ; e_{\alpha_1}+e_{\alpha_2} = e_{\beta_1}+e_{\beta_2}} E_{\alpha_1\alpha_2;\beta_1\beta_2}^*F_{\bmg\bmd} E_{\alpha_1\alpha_2;\beta_1 \beta_2} \ ,
$$
which for $\gamma_1\not=\gamma_2$  equals
$$
 \frac{\delta_{\gamma_1\delta_1} \delta_{\gamma_2\delta_2}}{2} 
 [f_{\gamma_1\gamma_1}\otimes f_{\gamma_2\gamma_2}  \otimes_{j=3}^Nf_{\gamma_j \delta_j} + f_{\gamma_2\gamma_2}\otimes f_{\gamma_1\gamma_1}  \otimes_{j=3}^Nf_{\gamma_j \delta_j}]
$$
$$
=\frac{\delta_{\gamma_1\delta_1} \delta_{\gamma_2\delta_2}}{2} 
 [f_{\gamma_1\delta_1}\otimes f_{\gamma_2\delta_2}  \otimes_{j=3}^Nf_{\gamma_j \delta_j} + f_{\gamma_2\delta_2}\otimes f_{\gamma_1\delta_1}  \otimes_{j=3}^Nf_{\gamma_j \delta_j}]
 $$
$$
=\frac{\delta_{\gamma_1\delta_1} \delta_{\gamma_2\delta_2}}{2} 
 [F_{\bmg\bmd} + S_{12}F_{\bmg\bmd}S_{12}]
 $$
 Hence we have
 \begin{eqnarray*}
 \cQ_{12}(F_{\bmg\bmd}) &=& \begin{cases}\delta_{\gamma_1\delta_1}\delta_{\gamma_2\delta_2} \delta_{\gamma_1\gamma_2}F_{\bmg\bmd}& \ {\rm if} \  \gamma_1 =\gamma_2 \\ \frac{\delta_{\gamma_1\delta_1} \delta_{\gamma_2\delta_2}}{2} 
 [F_{\bmg\bmd} + S_{12}F_{\bmg\bmd}S_{12}] & \ {\rm if} \ \gamma_1 \not= \gamma_2  \end{cases} \\ 
 &=& \frac{\delta_{\gamma_1\delta_1} \delta_{\gamma_2\delta_2}}{2} 
 [F_{\bmg\bmd} + S_{12}F_{\bmg\bmd}S_{12}]
 \end{eqnarray*}
since for $\gamma_1=\gamma_2$ the swap operator $S_{12}$ acts trivially. The same reasoning yields the analogous result for all other $m<n$. Summing and then using \eqref{quankac9} proves the lemma. 
\end{proof}

\begin{proof}[Proof of Theorem~\ref{GrLem2}] Fix $\bmb\in \G_M$, and  let $F_{\bmg,\bmd} = U_{S;\bmz\bmet}(g_{\bmb})$ so that $\bmg$ and $\bmd$ are given by \eqref{UnitDef}. Then $C_{\bmg,\bmd} = M$, and from  Lemma~\ref{WalkLem}
\begin{eqnarray*}
N^{-1} {\binom{N}{2}}\cL_N(F_{\bmg\bmd} ) 
&=& \frac 12\sum_{m<n} \delta_{\gamma_m\delta_m} \delta_{\gamma_n\delta_n}  [F_{\bmg\bmd} + F_{S_{mn}(\bmg)S_{mn}(\bmd)}]   -  F_{\bmg\bmd}\\
&=& \frac 12\sum_{m<n} \delta_{\gamma_m\delta_m} \delta_{\gamma_n\delta_n}  [ F_{S_{mn}(\bmg)S_{mn}(\bmd)} -F_{\bmg\bmd}]  \\
&-&  \left({\binom{N}{2}} - {\binom{M}{2}}\right) F_{\bmg\bmd}
\end{eqnarray*}
Now observe that
\begin{multline*}
U_{S;\bmz\bmet}\left(\frac 12 \sum_{m<n} \delta_{\gamma_m\delta_m} \delta_{\gamma_n\delta_n}  [ F_{S_{mn}(\bmg)S_{mn}(\bmd)} -F_{\bmg\bmd}] \right)U_{S;\bmz\bmet}^*(\bma)\\ = 
 \sum_{\{\bma'\ :\ [\bma,\bma']\in \E_M\}} g_\bmb(\bma') - v(\bma) g_\bmb(\bma)  = -\Delta_{\G_M}g_\bmb(\bma) \ ,
\end{multline*}
and
$$
{\binom{N}{2}} - {\binom{M}{2}} = \frac{(N^2- N)- (N^2-(2r+1)N +r(r+1)}{2} = rN - \frac{r(r+1)}{2}\ .
$$

\end{proof}

\section{The spectrum of the Graph Laplacian on the multislice}

The graph $\G_N$ is known as the {\em multislice} in analogy with the sliced Boolean cube $\{0,1\}^N$. The  
pair-transposition random walk on the vertices of the Boolean cube preserves the sum of the coordinates belonging to a vertex in  
$\{0,1\}^N$, and hence the paths of the walk stay in, and eventually cover, the ``slices'' of  
$\{0,1\}^N$ corresponding to the $N+1$ possible values of the sum of the coordinates.  

When $h$ is non-degenerate, we identify $\V_N = \{1,\dots,d\}^N$ with  $\{e_1,\dots,e_d\}^N$. Then due to energy conservation, our pair transposition walk on $\V_N$  conserves the value of the sum of the coordinates. The possible values of this sum are the eigenvalues of $H_N$,
and we have seen that these are indexed by the occupation vectors ${\bf k}(\bma)$, $\bma\in \G_N$.  Under the strong non-degeneracy condition, there is a one-to-one correspondence between eigenvalues $E$ of $H_N$,  and occupation vectors ${\bf k}  =(k_1,\dots,k_d)$ such that each $k_j$ is a non-negative integer and $\sum_{j=1}^dk_j =N$.  For each such ${\bf k}$, define $\G_{N,{\bf k}}$ to be the subgraph of $\G_N$ consisting of those vertices  $\bma\in \G_N$ for which ${\bf k}(\bma) = {\bf k}$. It is easy to see that these are precisely the connected components of $\G_N$. Thus $\Delta_{\G_N}$ is the direct sum of the operators $\Delta_{\G_N,{\bf k}}$.

\begin{theorem}\label{main} For all $N\geq 2 $, all $d\geq 2$, and all ${\bf k} = (k_1,\dots,k_{d}) \in \Z_{\geq 0}^d$ with 
\begin{equation}\label{full0}
\sum_{m=1}^{d}k_m = N\qquad {\rm and}\quad \max\{k_0,\dots,k_{r-1}\} < N\ ,
\end{equation}
so that $\G_{N,{\bf k}}$ is not trivial, the spectral gap $ \Gamma_{N,{\bf k}}$ of $\Delta_{\G_{N,{\bf k}}}$  is given by
\begin{equation*}
  \Gamma_{N,{\bf k}} = N\ .
\end{equation*}
\end{theorem} 

This together with our previous results yields
\begin{theorem}\label{Main2} 
For all $N\geq 2 $, all $d\geq 2$, and all ${\bf k} = (k_1,\dots,k_{d}) \in \Z_{\geq 0}^d$ with \eqref{full0} satisfied, the spectral gap of of
$\cL_{N,{\bf k}}$ is $N/(N-1)$. 
\end{theorem}

\begin{proof}
Now fix some $E\in \sigma(H_N)$ such that the corresponding occupation vector ${\bf k}$ satisfies the non-triviality condition of 
Theorem~\ref{main}.  Let $\H_{N,{\bf k}}$ be the corresponding eigenspace of $\H_N$, and let $P_{N,{\bf k}}$ be the 
orthogonal projection onto $\H_{N,{\bf k}}$. Then 
$\B(\H_{N,{\bf k}})$ is invariant under $\cL_N$, and
$\cL_N(P_{N,{\bf k}}) = 0$. Let $\cL_{N,{\bf k}}$ denote the restriction of $\cL_N$ to $\B(\H_{N,{\bf k}})$.

From Lemma~\ref{DecompLem}, we have
$$
\B(\H_{N,{\bf k}}) =  \bigoplus_{S\subset\{1,\dots,N\}\ ,\ (\bmz,\bmet)\in E_S}     \B_{S;\bmz,\bmet}
\ .
$$

For $S = \{1,\dots,N\}$, the eigenvalues  of $\cL_{N,{\bf k}}$ on $ \B_{S;\bmz,\bmet} $ are precisely $-1/(N-1)$ times the eigenvalues of
$\G_{N<{\bf k}}$, and hence the spectral gap in this sector is $N/(N-1)$. 

For $S$ a strict subset of  $\{1,\dots,N\}$, $|S| \leq N-2$, since when ${\bf k}(\bmg) = {\bf k}(\bmd)$, it is impossible for the coincidence set $S$ to have cardinality $N-1$. Hence in Lemma~\ref{GrLem2}, $r \geq 2$, and hence the least eigenvalue of  $-\cL_{N,{\bf k}}$ on $\B_{S;\bmz,\bmet} $ is at least
$$\left(2r\frac{N}{N-1} \textcolor{red}{-}  \frac{r(r+1)}{N-1}\right) \geq 4  - \frac{2}{N-1} \ .$$
\end{proof}

It is easy to write down a set of eigenfunctions of $\Delta_{\G_{N,{\bf k}}}$ that have eigenvalue $N$. This will be very useful in the proof that $N$ is in fact
 the spectral gap.

\begin{definition}\label{Kspace}  Let ${\mathcal K}_{N,{\bf k}}$ denote the set of real-valued functions $g$ on $\{e_1,\dots,e_{d}\}$ such that 
\begin{equation}\label{gapspace1}
\sum_{m=1}^{d} k_m g(e_m) = 0\ .
\end{equation}
\end{definition}

\begin{lemma}\label{Nex}   Let $g$ be any non-zero function in ${\mathcal K}_{N,{\bf k}}$ as specified in Definition~\ref{Kspace}.
Fix $1\leq \ell \leq N$ and define a function $f$ on $\V_{N,{\bf k}}$ by $f(x) = g(x_\ell)$. Then
\begin{equation}\label{espace1}
\Delta_{\G_{N,{\bf k}}} f(x) = Nf(x)\ .
\end{equation}
\end{lemma} 

\begin{proof}   Note that $f(x) - f(\pi_{i,j}x) = 0$ unless $i=\ell$ or $j=\ell$, hence we need only consider such pair permutations when computing 
$L_{\G_{N,{\bf k}}} f(x)$. 
For any $0 \leq m \leq r-1$, consider $x\in  \V_{N,{\bf k}}$ such that $x_\ell = e_m$.  For each $n\neq m$, there are $k_n$ pair permutations such that when applied to $x$ yield the value $e_n$ in the $\ell$th place. Therefore,

\begin{eqnarray*}
\Delta_{\G_{N,{\bf k}}} f(x)  &=& \sum_{n\neq m} k_n (g(e_m) - g(e_n)) \\
&=& (N-k_m) g(e_m)   - \sum_{n\neq m} k_ng(e_n)\\
&=&  (N-k_m) g(e_m) + k_mg(e_x) = Ng(e_m) = Nf(x)\ ,
\end{eqnarray*}
where in the last line we have used $\sum_{n=0}^{r-1}k_n g(e_n) = 0$. Since  $m$ is arbitrary, \eqref{espace1} is proved. 
\end{proof} 

\begin{remark}\label{null}  Fix any function $g$ on $\{e_1,\dots,e_{d}\}$ such that  $\sum_{n=1}^{d}k_n g(e_n) = 0$. Then the $N$ functions
$\{g(x_1), \dots, g(x_N)\}$  are not linearly independent since for any $x\in \V_{N,{\bf k}}$,
\begin{equation}\label{flat}
\sum_{\ell=1}^N g(x_\ell) = \sum_{n=1}^{d} k_n g(e_n) = 0\ .
\end{equation}
Now consider any $\{g_1,\dots,g_N\} \subset {\mathcal K}_{M,{\bf k}}$, not all zero,  and define the function
$$f(x) = 
\sum_{\ell =1}^N g_\ell(x_\ell)\ .
$$
By Lemma~\ref{Nex}, $L_{\G_{N,{\bf k}}} f = Nf$. However, we can express $f$ in a simpler way:  Since by \eqref{flat},  $\sum_{\ell=1}^N g_N(x_\ell) = 0$,
$$
f(x) = f(x) - \left(\sum_{\ell=1}^N g_N(x_\ell) \right) = \sum_{\ell =1} ^{N-1} h_\ell(x_\ell)
$$
where for $1 \leq \ell \leq N-1$, $h_\ell = g_\ell - g_N$.  
\end{remark}

\section{Spectral properties of the graph Laplacian on the multislice}

For each $N$, there is a natural partial order on the set of vectors ${\bf k} = (k_1,\dots,k_{d})$ induced by the ``coarsening operation'' of ``merging'' energy levels \cite{S20}. For any $d'> d \geq 2$, let $\phi:\{0,\dots,d'\} \to \{0,\dots,d\}$ be any surjection. For ${\bf k} = (k_0,\dots,k_{d'})$ with $\sum_{n=0}^{d'}k_n = N$, 
define 
$$\phi({\bf k})_m := \sum_{n\ :\ \phi(n)= m}k_n\quad{\rm and}\quad \phi({\bf k}) = (\phi({\bf k})_1,\dots,\phi({\bf k})_{d})\ .$$
Then $\phi$ induces a map, also denoted by $\phi$, from $\{e_1,\dots,e_{d'}\}$ into itself by
\begin{equation}\label{phiA}
\phi(e_n) = e_{\phi(n)} \ , \quad 1 \leq n \leq d'\ .
\end{equation}

We say that one multislice $\V_{N,{\bf k}}$ is {\em coarser} than another $\V_{N,{\bf k}'}$ if for some $2 \leq d < d' \leq N$, there is a surjection
$\phi:\{1,\dots,d'\} \to \{1,\dots,d\}$ such that $\phi({\bf k}') = {\bf k}$. 

In this case $\phi$ induces a map from $\V_{N,{\bf k}}$ onto $\V_{N,{\bf k}'}$ given by
\begin{equation}\label{phiB}
(\phi(x))_\ell := \phi(x_\ell)
\end{equation}
with $\phi(x_\ell)$ defined by \eqref{phiA}. 

The relevance of such coarsenings to spectral problems was pointed out and exploited in \cite{DF90,DS93}. In our setting, suppose that  for some 
$\phi$ as above we have ${\bf k}= \phi({\bf k}')$. It is easy to see that for any real valued function $f$  on $\V_{N,{\bf k}}$, 
\begin{equation}\label{phiC}
(\Delta_{\V_{N,{\bf k}}} f) \circ \phi  =  \Delta_{\V_{N,{\bf k}'}} (f \circ \phi)\ .
\end{equation}
Let $f$ be an eigenfunction of $\Delta_{\G_{N,{\bf k}}}$ with eigenvalue $\lambda$. 
Then
\begin{equation}\label{phiD}
\lambda f\circ \phi  = (\Delta_{\V_{N,{\bf k}}} f) \circ \phi  =  \Delta_{\V_{N,{\bf k}'}} (f \circ \phi)\ .
\end{equation}
Since the map $\phi$ defined in \eqref{phiB} is surjective, $f\circ \phi$ is not 
identically zero, and hence $\lambda$ is also an eigenvalue of $\Delta_{\V_{N,{\bf k}}}$.
Therefore, when the multislice $\V_{N,{\bf k}}$ is  coarser than the multislice  $\V_{N,{\bf k}'}$, the spectrum of $\Delta_{\G_{N,{\bf k}}}$ is 
contained in the spectrum of  $\Delta_{\G_{N,{\bf k}'}}$.   In particular,
\begin{equation}\label{phiE}
\Gamma_{N,{\bf k}'} \geq \Gamma_{N,{\bf k}}\ .
\end{equation}
Whenever $\G_{N,{\bf k}}$ is not trivial; i.e., has at least one edge,  it is known that $\Gamma_{N,{\bf k}}  =N$, independent of ${\bf k}$.  
This result can be found in \cite[Lemma 1]{S20}, in which the proof refers to the deep proof by Caputo, Liggett and Richthammer  
\cite{CLR} of a famous conjecture of Aldous \cite{AF}.  This fact about the gap also follows from \eqref{phiE} and some computations that follow. Because of  
\eqref{phiE}, it suffices to know $\Gamma_{N,{\bf k}}$ for a few special choices of ${\bf k}$.

Among the graphs considered here, some are absolutely trivial. For example, if  $k_m = N$ for some $0 \leq m \leq r-1$, then $\V_{N,{\bf k}}$ 
is a singleton, and  the edge set $\E_{N,{\bf k}}$  is empty. In this case, the graph Laplacian is $0$, and there is only this one eigenvalue, an hence no gap.

Somewhat less trivial is the case in which $k_{m_0} =N-1$ for some $m_0$. Then $k_m =1$ for one value of $m \neq m_0$, and 
$k_n = 0$ for all $n \neq m,m_0$.
In this case, every vertex is related to every  other by a pair transposition, and hence  the graph is a complete graph with 
$N$ vertices, and therefore the spectral gap is $N$.

For $r=2$, we have seen we might as well take $\{e_0,e_1\} = \{0,1\}$, and then $\{e_0,\dots,e_{r-1}\} = \{0,1\}^N$, the Boolean $N$-cube. 
Taking $\G = \{0,1\}^N$ with adjacency defined as above. The connected components $|G_{N(k_0,k_1)}$ are known as {\em Johnson Graphs}, and  the full spectrum of the Laplacian $L_{\G_{N,(k_0,k_1)}}$ is known \cite{BCN,Fil16}. In particular, it is known that  the spectral gap is always $N$ independent of ${\bf k} = (k_0,k_1)$ assuming that both $k_0$ and $k_1$ are non-zero.

Finally, consider the case in which $r=N$ and ${\bf k} = (1,\dots,1)$; i.e., $k_m =1$ for all $m$. Then evidently $\V_{N,(1,\dots,1)}$ has $N!$ 
vertices and may be identified with $S_N$, the symmetric group on $N$ letters.  The spectrum of the corresponding graph Laplacian 
$L_{\V_{N,(1,\dots,1)}}$ has been studied using methods from group representation theory by Diaconis and Shahshahani \cite{DS81}.  
Their results provide complete information on all of the eigenvalues, and this is essential for their applications. One of their results is that 
 for all $N$, the spectral gap is $N$.  For this alone, one does not need so much machinery, and a  simple proof is given in \cite[Theorem 5.1]{CCL2}. 

Now consider any ${\bf k}$ such that $\G_{N,{\bf k}}$ is non trivial.  Evidently ${\bf k}$ is coarser than
${\bf k}_1 := (1,\dots,1)$, and going in the opposite direction, reducing to just two energies, there is a ${\bf k}_0 = (k_0,k_1)$  such that 
\begin{equation}\label{comp}
\Gamma_{N,{\bf k}_0} \leq \Gamma_{N,{\bf k}} \leq \Gamma_{N,{\bf k}_1}\ .
\end{equation}
Since $\G_{N,{\bf k}_0}$ is a Johnson graph, $\Gamma_{N,{\bf k}_0} =  N$  \cite{BCN}. Then since  $\Gamma_{N,{\bf k}_1} =N$,
 $\Gamma_{N,{\bf k}} = N$. 
This result is  explicitly noted in \cite[Lemma 1]{S20}, with an argument based on the proof of Aldous' conjecture \cite{AF} by Caputo, Liggett and Richthammer 
\cite{CLR}.

We now give a simple proof of Theorem~\ref{main} and somewhat more using ideas developed in \cite{CCL2}, and applied there to the case of the sliced Boolean cube.

\subsection{The induction for the lower bound on the gap}

Define $\mu_{N,{\bf k}}$ to be
the uniform probability measure on $\V_{N,{\bf k}}$. 
The Dirichlet form for $\Delta_{\G_{N,{\bf k}}}$ on $L^2(\mu_{N,{\bf k}})$ is
\begin{equation}\label{unscaled}
\frac12 \sum_{x\in \V_{N,{\bf k}}} \sum_{i<j} (f(x) - f(\pi_{i,j}x))^2 \mu_{N,{\bf k}}\ .
\end{equation}
It greatly simplifies the induction we shall carry out if we  normalize  so that the sum over pairs becomes an average. Hence we work with the Dirichlet form of 
${{N}\choose {2}}^{-1}\Delta_{\G_{N,{\bf k}}}$:
\begin{equation}\label{unscaled2}
\frac12 {{N}\choose {2}}^{-1} \sum_{x\in \V_{N,{\bf k}}} \sum_{i<j} (f(x) - f(\pi_{i,j}x))^2 \mu_{N,{\bf k}}\ .
\end{equation}

The key to our induction is the identity
\begin{equation}\label{ave}
 {{N}\choose{2}}^{-1}\sum_{i < j} (f(\pi_{i,j} x) - f(x))^2 = \frac1N \sum_{\ell =1}^N \left(  {{N-1}\choose{2}}^{-1}\sum_{i < j, i,j \neq \ell}   (f(\pi_{i,j} x) - f(x))^2\right)\ ,
 \end{equation}
 On the right, we have an average of $N$ terms, each of which leaves  one coordinate, $x_\ell$, unchanged since transpositions $\pi_{i,j}$ in which either $i=\ell$ or $j=\ell$ are not included.  If one thinks in terms of  processes defined by the Dirichlet forms, this identity will relate the dynamics for $N$ particles to the dynamics for $N-1$ particles. 

First, we make one more adjustment to the $N$ particle dynamics. The Dirichlet form \eqref{unscaled2} is  associated to a continuous time Markov jump process on $\V_{N,{\bf k}}$ of the following description:  
A Poisson clock is running with expected times between ``rings'' equal to $2$. When a ``ring'' occurs, a pair $(i,j)$, $i < j$, is chosen 
uniformly at random, and the state jumps from vertex $x$ to vertex $\pi_{i,j}x$.  For any given $1 \leq \ell \leq N$, the number  
of pairs $i<j$ containing $\ell$ is $N-1$, and hence the fraction of the jumps that change the state of the $\ell$th particle is $2/N$. In order  to have
that all particles update with an expected time of order 1, independent of $N$, we therefore multiply the Dirichlet form in \eqref{unscaled2} 
by $N$, to obtain a family  of processes, indexed by $N$, in which the expected waiting times for updates of each particle are of order $1$, independent of $N$. This is physically motivated, but as we shall see, it is also convenient for the induction. 

\begin{definition}\label{dirdef}
Define the Dirichlet form 
 \begin{eqnarray}\label{dirformX}
 \mathcal{D}_{N,k}(f,f) &=&  \frac{N}{2} {{N}\choose{2}}^{-1}\sum_{x\in \V_{N,{\bf k}}}\sum_{i < j} (f(\pi_{i,j}x) - f(x))^2 \mu_{N,{\bf k}}(x)\nonumber\\
 &=&  \frac{1}{N-1}\sum_{x\in \V_{N,{\bf k}}}\sum_{i < j} (f(\pi_{i,j}x) - f(x) )^2 \mu_{N,{\bf k}}(x)
 \end{eqnarray}
 where the pair permutations $\pi_{i,j}$ acts on $x$ by swapping the $i$ and $j$th entries.   Also define $\widehat{\Gamma}_{N,{\bf K}}$ to be the spectral gap associated to this Dirichlet from. That is
 \begin{equation}\label{gapdef}
\widehat{\Gamma}_{N,{\bf k}} =  \inf\left \{    \mathcal{D}_{N,k}(f,f)   \ :\   \|f\|_{L^2(\mu_{N,{\bf k}})} =1\ , \ \langle f,1\rangle_ {L^2(\mu_{N,{\bf k}})}  = 0\ .   \right\}
 \end{equation} 
\end{definition}

\begin{remark} Comparing with \eqref{unscaled} which gives the Dirichlet form of the graph Laplacian $\Delta_{N,{\bf k}}$, we see that its gap, $\Gamma_{N,{\bf k}}$ and $\widehat{\Gamma}_{N,{\bf k}}$ are related by
\begin{equation}\label{gapprel}
\widehat{\Gamma}_{N,{\bf k}}   = \frac{2}{N-1} \Gamma_{N,{\bf k}}\ .
\end{equation}
\end{remark}

To make use of \eqref{ave}, we first consider a graph $\G_{N,{\bf k}}$ where ${\bf k} = (k_1,\dots,k_{d})$  is such that 
\begin{equation}\label{kcond}
k_m \geq 1 \quad{\rm  for\ each} \quad 1 \leq m \leq d\ .
\end{equation}
Then $\G_{N,{\bf k}}$ is  a graph for $N$ particles that truly have $d$ different energy levels. If it were the case that $k_m= 0$ for some $m$, 
the energy $e_m$ would play no role, and the graph would be identical to  another graph with a reduced set of $r' < r$ energy levels and a net ${\bf k}'$
  such that  for all $0 \leq m \leq r'-1$, $k'_m\geq 1$.   Evidently $\G_{N,{\bf k}}$ and $\G_{N,{\bf k}'}$ are isomorphic and in 
  particular  have the same spectral gap.

 Now considering $\G_{N,{\bf k}}$ such that 
 \eqref{kcond} is satisfied, we specify a bijection of $\V_{N,{\bf k}}$ with a union of vertex sets of graphs for $N-1$ particles:
    For $0 \leq m \leq r-1$, define ${\bf k}^{(m)}$ to be obtained from ${\bf k}$ by replacing $k_m$ with $k_m-1$.  For each $1 \leq \ell \leq N$ we define a map
 \begin{equation}\label{GJ2}
 T_\ell : \left(\bigcup_{m=0}^{r-1} \V_{N-1,{\bf k}^{(m)}}\right) \to \V_{N,{\bf k}}
 \end{equation}
 by
 \begin{equation}\label{GJ3}
 T_\ell(x) = (x_1, \dots, x_{\ell-1}, e_m,x_\ell,\dots, x_{N-1})\quad{\rm for}\quad x\in  \V_{N-1,{\bf k}^{(m)}}
  \end{equation}
  with the obvious modifications for $\ell =1$ or $\ell = N-1$.   In other words, for $x\in \V_{N-1,{\bf k}^{(m)}}$, one simply inserts $e_m$ in the $\ell$th place, keeping the order of the remaining entries unchanged, and thus obtains an element of $\V_{N,{\bf k}}$. It is evident that this does indeed yield a bijection. 
  The set $\{x\in \V_{N,{\bf k}} \ :\ x_\ell = e_m\}$ is precisely the image of $\V_{N-1,{\bf k}^{(m)}}$ under $T_\ell$. 
  We shall prove
    
  \begin{theorem}\label{maxmain}  Let $\{e_1,\dots, e_{d}\}$ be given along with ${\bf k} = (k_0,\dots,k_{r-1})$ where each $k_m$ is a strictly positive integer and $\sum_{m=1}k_m = N$. Then the spectral gap $\Delta_{N,{\bf k}}$ for the corresponding Dirichlet form  specified in \eqref{dirformX} satisfies
 \begin{equation}\label{indlong}
 \widehat{\Gamma}_{N,{\bf k}} \geq     \frac{N(N-2)}{(N-1)^2} \min\{  \widehat{\Gamma}_{N-1,{\bf k}^{(m)}}  \ :\  0 \leq m \leq r-1 \}
 \end{equation}
  \end{theorem}

  We shall prove this to be a consequence of several lemmas, and begin by  explaining how we make use of \eqref{ave}

For each $1 \leq \ell \leq N$, and each $m\in \{0,\dots,r-1\}$, define
 \begin{eqnarray}\label{dirform2}
  \mathcal{D}^{\ell,m} _{N,{\bf k}}(f,f)
   &=& \frac{N-1}{2}{{N-1}\choose{2}}^{-1} \sum_{x\in \V_{N,{\bf k}}, x_\ell = e_m}\left(\sum_{i < j, i,j \neq \ell} (f(\pi_{i,j} x) - f(x))^2 \right) 
 \mu_{N-1,{\bf k}^{(m)}}\nonumber\\
 &=& \frac{1}{N-2} \sum_{x\in \V_{N,{\bf k}}, x_\ell = e_m}\left(\sum_{i < j, i,j \neq \ell} (f(\pi_{i,j} x) - f(x))^2 \right) \mu_{N-1,{\bf k}^{(m)}}\ .
 \end{eqnarray}
 
 Next, for each $1 \leq \ell \leq N$ define the operator $P_\ell$ on $L^2(\V_{N,{\bf k}})$ as follows: On the set $\{x\ :\ \ x_\ell = e_m\}$,
   $$ 
   P_\ell f(x) :=   \mu_{N,{\bf k}^{(m)}} \sum_{y \in \V_{N,{\bf k}}\ :\  y_\ell = e_m} f(y)   \ .
   $$
   The operator $P_\ell$ is just the orthogonal projection in $L^2(\mu_{N,{\bf k}})$ onto the subspace of functions that depend only on $x_\ell$. 
Therefore,  for each $m\in \{0,\dots,r-1\}$,
\begin{equation}\label{shift}
  \mathcal{D}^{\ell,m} _{N,{\bf k}}(f,f) = \mathcal{D}^{\ell,m} _{N,{\bf k}}(f - P_\ell f ,f - P_\ell f) \ .
 \end{equation}
 Also, note that
 \begin{equation}\label{GJG21}
 \mu_{N,{\bf k}} = \sum_{m=0}^{r-1} \frac{k_m}{N} \mu_{N-1,{\bf k}^{(m)}}\ . 
 \end{equation}
 Therefore, using the key identity \eqref{ave}, together with \eqref{shift} and \eqref{GJG21},
  \begin{equation}\label{GJG22}
  \mathcal{D}_{N,{\bf k}}(f,f)  = \frac1N \sum_{\ell=1}^N  \frac{N}{N-1}\sum_{m=0}^{r-1}   \mathcal{D}^{\ell,m} _{N,{\bf k}}(f-P_\ell f ,f- P_\ell f) \frac{k_m}{N}
   \end{equation}
 
Recall that  $ \widehat{\Gamma}_{N,{\bf k}}$ is  the spectral gap of the Dirichlet form defined in \eqref{dirformX}.   Since for each $m=0,\dots,r-1$, $f - P_\ell f$ is constant on the image of each $\V_{N-1,{\bf k}^{(m)}}$ under $T_\ell$, 
 $$
  \mathcal{D}^{\ell,m} _{N,{\bf k}}(f-P_\ell f ,f- P_\ell f) \geq \Delta_{N-1,{\bf k}^{(m)}} \|(f - P_\ell f)\circ T_\ell\|^2_{L^2(\mu_{N-1,{\bf k}^{(m)}})}\ .
 $$
 By \eqref{GJG21}, and the fact that each $P_\ell$ is an orthogonal projection,
 \begin{eqnarray*}
 \sum_{m=1}^{r-1}  \|(f - P_\ell f)\circ T_\ell\|^2_{L^2(\mu_{N-1,{\bf k}^{(m)}})}\frac{k_m}{N} &=& \|f -  P_\ell f\|_{L^2(\mu_{N,{\bf k}})} \\
&=&  \|f\|^2_{L^2(\mu_{N,{\bf k}})}   - \langle f, P_\ell  f \rangle_{L^2(\mu_{N,{\bf k}})}\ .
 \end{eqnarray*}
 Then taking $f$ to be a  normalized gap eigenfunction for $\Delta_{N,{\bf k}}$, we have from \eqref{GJG22} that
 \begin{equation}\label{PinductX}
   \widehat{\Gamma}_{N,{\bf k}} \geq \min\left\{  \widehat{\Gamma}_{N-1,{\bf k}^{(m)}} \ :\  0 \leq m \leq r-1\right\} \frac{N}{N-1}(1-\langle f, P f\rangle_{L^2(\mu_{N,{\bf k}})} )
\end{equation}
 where
 \begin{equation}\label{pop}
 P =\frac1N  \sum_{\ell=1}^N P_\ell\ .
 \end{equation}
 Note that $P$ is an average of orthogonal projections, and hence its spectrum lies in $[0,1]$. Any function $f$ that is an eigenvalue of 
 $P$ with the eigenvalue $1$ must be in the range of each of the projections $P_\ell$. However, the range of $P_\ell$ consists of functions 
 $f$ that depend on $x$ only through $x_\ell$. 
 For $N\geq 3$, the only functions $f$ on $\V_{N,{\bf k}}$ that have this property for all $\ell$ are the constant functions since for every 
 $x,y\in \V_{N,{\bf k}}$, there is a sequence of pair transpositions that takes $x$ to $y$, and since $N\geq 3$, each such transposition 
 leaves one coordinate unchanged, and hence leaves the value of $f$ unchanged. Hence $f(x) = f(y)$. 
Therefore, for $N\geq 3$,  $1$ is an eigenvalue of $P$ of multiplicity one with the eigenspace being the constant functions. 
 Since the gap eigenfunction $f$ is orthogonal to the constants in 
 $L^2(\mu_{N,{\bf k}})$, we must have $\langle f, P f\rangle < 1$. In fact, this quantity can be no larger that the next largest eigenvalue of $P$:
 
 \begin{definition} Let $N\geq 3$ and let $\lambda_{N,{\bf k}}$ denote the second largest eigenvalue  of $P$
 \begin{equation}\label{Pgapdef}
 \lambda_{N,{\bf k}}   = \sup\left\{  \langle h,P h\rangle_{L^2(\mu_{N,{\bf k}})}  \ :\  \|h\|_{L^2(\mu_{N,{\bf k}})}  =1\ ,\ \langle h,1\rangle_{L^2(\mu_{N,{\bf k}})} = 0 \right\}\ .
 \end{equation}
 \end{definition} 
 Therefore, \eqref{PinductX} becomes
 \begin{equation}\label{PinductXL}
   \widehat{\Gamma}_{N,{\bf k}} \geq \min\left\{  \widehat{\Gamma}_{N-1,{\bf k}^{(m)}} \ :\  0 \leq m \leq r-1\right\} \frac{N}{N-1}(1-\lambda_{N,{\bf k}})\ .
\end{equation}

The next lemma renders \eqref{PinductXL} completely explicit and yields the proof of Theorem~\ref{maxmain}.

  \begin{lemma}\label{Plm}  Let $N\geq 3$.  The spectrum of $P$ is the set $\{0, (N-1)^{-1}, 1\}$. In particular,
  \begin{equation}\label{lam}
  \lambda_{N,{\bf k}} = \frac{1}{N-1}\ .
  \end{equation}
  Moreover, the eigenspace corresponding to $1$ consists of the constant functions, and 
  the eigenspace corresponding to $(N-1)^{-1}$   has dimension $(r-1)(N-1)$ and a basis for it   is the set  of   functions of the form
  \begin{equation}\label{Plm1}
  f_{m,\ell}(x) = g_m(x_\ell)  \quad 1 \leq m \leq r-1 \quad {\rm and}\quad 1 \leq \ell  \leq N-1\ 
  \end{equation}
  where $\{g_1,\dots,g_{r-1}\}$ is a basis for  $\mathcal{K}_{N,{\bf k}}$. 
  \end{lemma}
  
 The proof Lemma~\ref{Plm} that we give is patterned on the proof of Lemma~2.16 of \cite{CCL3}.  It involves a simpler operator $K$ that, like $P$, is constructed out of the projections $P_\ell$.  K is an operator on functions of a single variable in $\{e_1,\dots,e_{d}\}$:  Let  $\nu_{N,{\bf k}}$ denote the probability measure on  $\{e_1,\dots,e_{d}\}$
 that is the image of the uniform probability measure on $\V_N$ under the map $\pi_N:(x_1,\dots,x_N)\mapsto x_N$. It is easy to see that 
  \begin{equation}\label{GJ4}
   \nu_{N,{\bf k}}(\{e_m\}) = \frac{k_m}{N}\ .
   \end{equation} 
Define the operator $K$ on $L^2(\nu_{N,{\bf k}})$ by $Kh = P_1h\circ \pi_N$. Note that for any function $h$ on $\{e_1,\dots,e_{d}\}$, $h\circ \pi_N$ is a function on $\V_N$, depending only on the $N$th coordinate, and then $P_1h\circ \pi_N$, which is a function of the first coordinate only, may be identified with a function on $\{e_1,\dots,e_{d}\}$, and this function is, by definition, $Kh$.  There is nothing special about $1$ and $N$, and pairs of distinct indices yields the same operator by symmetry. 

By the definition of $\nu_{N,{\bf k}}$ in terms of $\mu_{N,{\bf k}}$, we have the following formula for $K$ which shows that it is self-adjoint on  $L^2(\nu_{N,{\bf k}})$:
  \begin{equation}\label{Kdef1X}
\langle g, K h\rangle_{L^2(\nu_{N,{\bf k}})} :=  \sum_{x\in \V_{N,{\bf k}}} g(x_1)h(x_N) \mu_{N,{\bf k}}(x)\ .
 \end{equation}

 \begin{lemma}\label{Klm}  The spectrum of the  $K$ on $L^2(\nu_{N,{\bf k}})$ is $\{1, -1/(N-1)\}$. The eigenspace corresponding to the eigenvalue $1$ consists of the constant functions on $\{e_0,\dots,e_{r-1}\}$, and the eigenspace corresponding to the eigenvalue $-1/(N-1)$ is the space
  consisting of functions $g$ on $\{e_0,\dots,e_{r-1}\}$ such that $\sum_{m=1}^{r-1}k_mg(e_m) = 0$. 
 \end{lemma}
 
\begin{proof}
Since since for $m\neq n$, there are
${\displaystyle\frac{(N-2)!}{k_0! \cdots k_{r-1}!}k_mk_n}$
vertices with $x_1 = e_m$ and $x_N = e_n$, and there are
${\displaystyle \frac{(N-2)!}{k_0! \cdots k_{r-1}!}k_m(k_m-1)}$
vertices with $x_1 = x_N = e_m$, working from the right side of  \eqref{Kdef1X}, we find
 \begin{eqnarray*}
 \sum_{x\in \V_{N,{\bf k}}} g(x_1)h(x_N) \mu_{N,{\bf k}}(x) &=& \frac{1}{N(N-1)} \sum_{m=0}^{r-1}g(e_m)h(e_m)k_m(k_m-1) \\
 &+& \frac{1}{N(N-1)} \sum_{m \neq n} g(e_m)h(e_n) k_mk_n \ .
 \end{eqnarray*}
 From here it follows easily that
 $Kh(e_m) = \sum_{n=0}^{r-1} K_{m,n} h(e_n)$
 where
 \begin{equation}\label{KmatrixX2}
(N-1) K_{m,n} = \begin{cases}  k_n-1 & n=m\\ k_n & n\neq m\end{cases}\ .
 \end{equation}
    Therefore,  $(N-1)K   = -1 + L$ where $L$ is the rank one matrix each of whose rows is $(k_0,\dots,k_{r-1})$.   
    Evidently, the non-zero constant vectors are eigenvectors of $K$  with eigenvalue $1$, and every non-zero vector 
    orthogonal to the constants in $L^2(\nu_{N,{\bf k}})$ is an eigenvector with eigenvalue $-1/(N-1)$. 
 \end{proof}

 Note that the space ${\mathcal K}_{N,{\bf k}}$ that  figures in Theorem~\ref{main} is precisely   the eigenspace of $K$ corresponding to the eigenvalue $-1/(N-1)$.

 We have seen in Remark~\ref{null}  that for $\{g_1,\dots,g_N\} \subset {\mathcal K}_{N,{\bf k}}$, $\{g_1(x_1),\dots,g_N(x_N)\}$ 
 need not be linearly independent even if each $g_\ell$ is non-zero.   The following lemma shows that the only way linear 
 independence can fail is the way indicated in Remark~\ref{null}
 
\begin{lemma}\label{nulllm}  Let  $\{g_1,\dots,g_N\} \subset {\mathcal K}_{N,{\bf k}}$ Then
$\sum_{\ell=1}^N  g_\ell(x_\ell) = 0$ for all $x\in \V_{N,{\bf k}}$ if and only if all of the functions $g_1,\dots,g_N$ are the same. 
\end{lemma} 

\begin{proof}  By what is explained in Remark~\ref{null}, it suffices to show that if $\sum_{\ell=1}^N  g_\ell(x_\ell) = 0$ for all $x\in \V_{N,{\bf k}}$, then all of the 
functions $g_1,\dots,g_N$ are the same. Hence we assume $\sum_{\ell=1}^N  g_\ell(x_\ell) = 0$ and compute
\begin{eqnarray*}
0 &=& \left\Vert \sum_{\ell=1}^N  g_\ell(x_\ell) \right\Vert^2_{L^2(\mu_{N,{\bf k}})}  
= \sum_{i=1}^N \left( \|g_i\|^2_{L^2(\nu_{N,{\bf k}})}  +  \sum_{\ell \neq i} \langle g_i, K g_\ell\rangle_{L^2(\nu_{N,{\bf k}})}\right)\\
&=&\sum_{i=1}^N \left( \|g_i\|^2_{L^2(\nu_{N,{\bf k}})} - \frac{1}{N-1}   \sum_{\ell \neq i} \langle g_i, g_\ell\rangle_{L^2(\nu_{N,{\bf k}})}\right)\\
&\geq&  \sum_{i=1}^N \left( \|g_i\|^2_{L^2(\nu_{N,{\bf k}})}  -  \frac{1}{N-1} \sum_{\ell \neq i} \| g_i\|_{L^2(\nu_{N,{\bf k}})}  \|g_\ell\|_{L^2(\nu_{N,{\bf k}})}\right)\\
&\geq&  \sum_{i=1}^N \left( \|g_i\|^2_{L^2(\nu_{N,{\bf k}})}  -  \frac{1}{N-1} \sum_{\ell \neq i} \frac12(\| g_i\|^2_{L^2(\nu_{N,{\bf k}})} +
  \|g_\ell\|^2_{L^2(\nu_{N,{\bf k}})})\right) = 0\ .
\end{eqnarray*}
Hence both inequalities must be equalities. The first inequality is the Schwarz inequality and a worst case assumption on the signs, 
and equality here entails that the  functions $g_1,\dots,g_N$ are  are all positive multiples of one another. The second inequality is the 
arithmetic-geometric mean inequality, and equality here entails that $\|g_i\|_{L^2(\nu_{N,{\bf k}})}  = \|g_\ell\|_{L^2(\nu_{N,{\bf k}})} $ for all $i$ and $\ell$. 
\end{proof} 
      
%
%

  \begin{proof}[Proof of Lemma~\ref{Plm}]  Suppose that $f$ is an eigenfunction of $P$ with eigenvalue $\lambda>0$.   Write
  $$ P_\ell f(x) =: h_\ell(x_\ell) = h_\ell\circ\pi_\ell(x)\ ,$$
  where $\pi_\ell(x) = x_\ell$.
  Define 
  $\vec h := (h_1,\dots,h_N)$
  viewed as an element of $\C^N\otimes  L^2(\nu_{N,{\bf k}})$.  Since $f$ is an eigenfunction of $P$,  $f\neq 0$ and 
  $N\lambda \sum_{\ell=1}^N h_\ell(x_\ell) = f(x)$.  Therefore  $\vec h \neq 0$, and the map $f\mapsto \vec h$ is injective from  
  the eigenspace of $P$ corresponding to eigenvalue $\lambda>0$  into 
  $\C^N\otimes  L^2(\nu_{N,{\bf k}})$.  Since this is true for every eigenvalue $\lambda>0$ of $P$, the map $f\mapsto \vec h$ is injective from 
  ${\rm ker}(P)^\perp$ into $\C^N\otimes  L^2(\nu_{N,{\bf k}})$.

  Applying $P_i$ to both sides of $\lambda f = Pf$,  and using the formula $Kh_\ell = P_i h_\ell\circ \pi_\ell$,
    \begin{equation}\label{fpr1}
 (N \lambda -1)  h_i  =  \sum_{\ell \neq i} K h_\ell \ .
  \end{equation}
  Define $M$  to be the $N\times N$ matrix with 
 ${\displaystyle M_{i,\ell} = \begin{cases} 0 & i = \ell\\ 1 & i\neq \ell
 \end{cases}}$.
 Then since \eqref{fpr1} is valid for each $i$,  we have
 $$
  (N \lambda -1) \vec h = M\otimes K \vec h\ .
 $$
 Since $\vec h\neq 0$, $(N \lambda -1)$ is an eigenvalue of $M\otimes K$. Evidently the spectrum of $M$ is $\{-1,N-1\}$, and by Lemma~\ref{Klm}, 
 the spectrum of $K$ is $\{1,-1/(N-1)\}$.  Then the spectrum of $M\otimes K$ is $\{ -1, 1/(N-1), N-1\}$.  It follows that if  $\lambda> 0$, then 
  either $\lambda = 1$ or $\lambda = 1/(N-1)$.  Since $(N\lambda -1)\in \{ -1, 1/(N-1), N-1\}$ is equivalent to
  $$\lambda \in \{0, 1/(N-1), 1\}\ ,$$
  this proves that the only eigenvalue $\lambda$ of $P$ with $0 <\lambda <1$ is $\lambda = 1/(N-1)$, and hence \eqref{lam} is proved. 
  
  Moreover, we have seen that the map $f\mapsto \vec h = (h_1,\dots, h_N)$ is injective from the $1/(N-1)$ eigenspace of $P$ into the $1/(N-1)$eigenspace of 
  $M\otimes K$, and this eigenspace is the product of the $-1$ eigenspace of $M$ and the $-1/(N-1)$ eigenspace of $K$.  
  Evidently, the dimension of this eigenspace is $(N-1)(r-1)$.  Hence the dimension of the $1/(N-1)$ eigenspace of 
  $P$ can not have a dimension any higher than this. 
  
To see that it is not any lower, consider any nonzero $g\in {\mathcal K}_{N,{\bf k}}$, and any $1 \leq \ell \leq N$, and let $f(x) = g(x_\ell)$. We then compute
\begin{eqnarray*}
Pf(x) &=& \frac1N\left(g(x_\ell) + \sum_{i\neq \ell} Kg(x_i)\right)   = \frac1N\left( g(x_\ell) -\frac{1}{N-1}  \sum_{i\neq \ell }g(x_i) \right)\\
&=& \frac1N\left( \frac{N}{N-1} g(x_\ell) - \frac{1}{N-1}\sum_{i=1}^N g(x_i)\right) = \frac{1}{N-1} g(x_\ell) = \frac{1}{N-1} f(x)\ .
\end{eqnarray*}

 Hence all such functions belong to the eigenspace corresponding to the eigenvalue $1/(N-1)$. Choose a basis $\{g_1,\dots,g_{r-1}\}$
 of $\mathcal{K}_{N,{\bf k}}$.   Let $f$ be a nontrivial linear combination of the $(N-1)(r-1)$ functions 
 $f_{m,\ell} = g_m(x_\ell)$, $1 \leq m \leq r-1$ and $1 \leq \ell \leq N-1$. The result is a function of the form 
 $\sum_{\ell =0} h_\ell(x_\ell)$ with each $h_\ell \in \mathcal{K}_{N,{\bf k}}$, at least one of which is 
 non-zero. 
 Then by Lemma~\ref{nulllm}, this cannot vanish identically, and hence the specified set of $(N-1)(r-1)$ functions is
  linearly independent. Hence they constitute a basis for the eigenspace corresponding to the eigenvalue $1/(N-1)$.
eigenspace corresponding to the eigenvalue $1/(N-1)$
\end{proof}
 
 \begin{proof}[Proof of Theorem~\ref{maxmain}]  This now follows directly from \eqref{PinductXL} and \eqref{lam}. 
 \end{proof}

\subsection{Proof of  Theorem~\ref{main}}

To use our inductive relation \eqref{PinductXL}, we need to know the values of $\Gamma_{N,{\bf k}}$ for small $N$. For some values of ${\bf k}$, this is trivial even for large $N$:   If $\max\{ k_m\ :\ 0 \leq m \leq r-1 \} = N$, $\G_{N,{\bf k}}$ has a single vertex and no edges. We exclude these trivial cases, and going forward suppose that
\begin{equation}\label{dirformY1}
\max\{ k_m\ :\ 0 \leq m \leq r-1 \} \leq  N-1 \ .
\end{equation}
If there is equality in \eqref{dirformY1}, then every vertex  $x= (x_1,\dots,x_N)$ in  $\V_{N,{\bf k}}$ has all but one of the entries $x_j$ the same, 
and exactly one that is not. All of these are related to one another by a pair transposition, and hence in this case, $\G_{N,{\bf k}}$ is a 
complete graph with $N$ vertices, and hence 
\begin{equation}\label{dirformY12}
 \Gamma_{N,{\bf k}}  = N \qquad{\rm and}\qquad  \widehat{\Gamma}_{N,{\bf k}}  =\frac{2N}{N-1}  
\end{equation}
for all such $k$. 

Finally, as we have observed before, it suffices to consider graphs $\G_{N,{\bf k}}$ for which 
$k_m \geq 1$ for  $0 \leq m \leq r-1$, since otherwise the graph is the same as one with a smaller set of energies that does satisfy such a condition. 

Consider $N=2$.  The only non-trivial choice for ${\bf k}$  is with $r=2$ and  ${\bf k} = (1,1)$. There are two vertices 
$(e_0,e_1)$ and $(e_1,e_0)$ and the single edge connects them. This is a complete graph, and hence $\Gamma_{2,(1,1)} = 2$  
and $\widehat{\Gamma}_{2,(1,1)} = 4$.  In summary, for $N=2$, the is only one non-trivial choice of ${\bf k}$, and for this choice, $\widehat{\Gamma}_{2,{\bf k}} = 4$. 

Next consider $N=3$.  The non-trivial choices for ${\bf k}$ are, with $r=2$, ${\bf k} = (1,2)$ and ${\bf k} = (2,1)$, both of which 
are complete graphs, and for $r=3$, ${\bf k} = (1,1,1)$. By Theorem~\ref{maxmain}
$$
\widehat{\Gamma}_{3,(1,1,1)} \geq \frac34 \min\{ \widehat{\Gamma}_{2,(0,1,1)}\ , \widehat{\Gamma}_{2,(1,0,1)}\ , \ \widehat{\Gamma}_{2,(1,1,0)}\ \}
$$
But evidently $\widehat{\Gamma}_{2,(0,1,1)}=\widehat{\Gamma}_{2,(1,0,1)} =  \widehat{\Gamma}_{2,(1,1,0)} =  \widehat{\Gamma}_{2,(1,1)} = 4$.
Therefore
$$
\widehat{\Gamma}_{3,(1,1,1)} \geq \frac{3}{4} 4 = 3\ .
$$ 
Since $\G_{3,(2,1)}$ and $\G_{3,(1,2)}$ are complete, $\Gamma_{3,(2,1)} = \Gamma_{3,(1/2)} =3$, and since  $2/(N-1) =1$ for 
$N=3$, we also have
$\widehat{\Gamma}_{3,(2,1)} = \widehat{\Gamma}_{3,(1,2)} =3$. In summary, for $N=3$ and all non-trivial choices of ${\bf k}$, $\widehat{\Gamma}_{3,{\bf k}} =3$. 

\begin{proof}[Proof of Theorem~\ref{main}]

Let $N\geq 3$ be an integer. We make the inductive hypothesis that for $M =N-1$, $\widehat{\Gamma}_{M,{\bf k}} = 2M/(M-1)$ for all ${\bf k}$ such that 
$\G_{M,{\bf k}}$ is non-trivial.   By the remarks made above, this is valid for $M=2$ and $M=3$. Now consider ${\bf k}$ such that 
$\G_{N,{\bf k}}$ is non-trivial and such that $k_m \geq 1$ for $0 \leq m \leq r-1$ which, as explained above, we may assume 
without loss of generality. Then, by Theorem~\ref{maxmain},
\eqref{PinductXL} and the inductive hypothesis yield
$$
\widehat{\Gamma}_{N,{\bf k}} \geq  \frac{N(N-2)}{(N-1)^2}  \frac{2(N-1)}{N-2} = \frac{2N}{N-1}\ . 
$$
However, by Lemma~\ref{Nex} and \eqref{gapprel}, 
$$\widehat{\Gamma}_{N,{\bf k}} \leq  \frac{2N}{N-1}\quad {\rm and\ hence}\quad  \widehat{\Gamma}_{N,{\bf k}} =  \frac{2N}{N-1}\ $$
for all $N$ and all ${\bf k}$ such that $\G_{N,{\bf k}}$ is non-trivial.   By \eqref{gapprel} once more,  this proves that $\Gamma_{N,{\bf k}} = N$.  
for all $N$ and all ${\bf k}$ such that $\G_{N,{\bf k}}$ is non-trivial.

Now take $f$ to a  normalized gap eigenfunction for $\widehat{\Gamma}_{N,{\bf k}}$.   We have from \eqref{PinductX} that if $f$ is not a gap 
eigenfunction of $P$, there there is  strict inequality in \eqref{PinductXL}, and this in turn would yield
$$
\widehat{\Gamma}_{N,{\bf k}} >  \frac{N(N-2)}{(N-1)^2}  \frac{2(N-1)}{N-2} = \frac{2N}{N-1}\ . 
$$
This contradiction shows that every gap eigenfunction for $\cL_{N,{\bf k}}$ is a gap eigenfunction of $P$.   However, 
Lemma~\ref{Plm} provides a complete description of the gap eigenspace of $P$, and Lemma~\ref{Nex} shows that every 
gap eigenfunction of $P$ is a gap eigenfunction of $\cL_{N,{\bf k}}$. 
\end{proof}

\section{Relative entropy dissipation}

Let $\tau_{N,{\bf k}}$ denote  the normalized trace on $\H_{N,{\bf k}}$. That is, for $X\in \B(\H_{N,{\bf k}})$,
$$
\tau_{N,{\bf k}}(X) := \frac{1}{{\rm dim}(\H_{N,{\bf k}})}\tr[X]\ .
$$
We say that $\varrho\in \B(\H_{N,{\bf k}})$ is a normalized density matrix in case $\varrho \geq 0$ and $\tau_{N,{\bf k}}(\varrho) =1$. 
Every quantum state on $\B(\H_{N,{\bf k}})$ has a representation of the form  $X\mapsto  \tau_{N,{\bf k}}(X\varrho)$ for some uniquely determined normalized density matrix $\rho$.

Recall that  $P_{N,{\bf k}}$ is the  orthogonal projection onto $\H_{N,{\bf k}}$.  Then $P_{N,{\bf k}}$ is a normalized density matrix in 
$\B(\H_{N,{\bf k}})$, and as we have seen 
$\cL_{N,{\bf k}}(P_{N,{\bf k}}) =0$. That is $P_{N,{\bf k}}$ is the equilibrium state in the sector $\B(\H_{N,{\bf k}})$.

The relative entropy of $\varrho$ with respect to the equilibrium state $P_{N,{\bf k}}$ is the  quantity
$$
D(\varrho||P_{N,{\bf k}}) := \tau_{N,{\bf k}}( \varrho(\log \varrho - \log P_{N,{\bf k}})) = \tau_{N,{\bf k}} (\varrho \log \varrho)\ .
$$

The are two entropy  inequalities that  are very useful for studying the approach to equilibrium, each of which implies a spectral gap inequality.  
The first is a generalized logarithmic Sobolev inequality and it takes the following form in our context:  Let $C_{N,{\bf k}}$ be the 
such that
\begin{equation}\label{GLS}
\tr[\varrho \log \varrho] \leq  -C_{N,{\bf k}} \tr[ \log \varrho \cL_{N,{\bf k}} \varrho ]
\end{equation}
for all normalized density matrices  $\varrho$  in $\B(\H_{N,{\bf k}})$. 

In this case, we have that for any normalized density matrix $\varrho_0$,
$$
D(e^{t{\cL_{N,{\bf k}}}}\varrho_0||P_{N,{\bf k}})   \leq e^{-t C_{N,{\bf k}} }D(\varrho_0||P_{N,{\bf k}}) \ .
$$
The inequality \eqref{GLS} is known as a {\em generalized logarithmic Sobolev inequality}. 

The existence of a finite constant 
$C_{N,{\bf k}}$ such that \eqref{GLS} holds is trivial, but determining the dependence on ${\bf k}$ and $N$ is not.  
In the classical sector,  or, what is the same thing, for the corresponding walk on slices of the multislice, it is known that that 
$1/2 \leq C_{N,{\bf k}} \leq 1$ for all non-trivial ${\bf k}$ and $N$. See \cite[Section 3.4]{S20} where this is deduced from a 
comparison argument and a result of Caputo, Dai Pra and Posta \cite{CDP}.  It is natural to conjecture that a similar result is valid for the quantum model, and this is the subject of current research.


\begin{thebibliography}{30}

\small{

 
\bibitem{AAKV} D.~Aharonov, A.~Ambainis, J.~Kempe and U.~Vazirani, \textit{Quantum walks on graphs}, 
STOC '01: Proceedings of the thirty-third annual ACM symposium on Theory of computing, (2001) 50--59


\bibitem{AF} D.~Aldous and  J.~Fill, \textit{Reversible Markov Chains and Random Walks on Graphs} 
http://www.stat.berkeley.edu/$\tilde{\phantom{.}}$aldous/RWG/book.html


\bibitem{BCN} A.~E.~Brower, A.~M.~Cohen and A. Neumaier, \textit{Distance-Regular Graphs}, Springer Berlin Heidelberg, 1989. 

\bibitem{CDP}  P.~Caputo, P.~Dai Pra, and G.~Posta, \textit{Convex entropy decay via the Bochner-Bakry-Emery approach}, Ann. Inst. Henri Poincar\'e, Probab. Stat. {\bf 45} (2009),  734--753.

\bibitem{CLR} P.~Caputo, T.~Liggett and T.~Richthammer, \textit{Proof of Aldous' spectral gap conjecture}, J. Amer. Math. Soc. {\bf 23} (2010), 831--851.




\bibitem{CCL2} E.~A.~Carlen, M.~C.~Carvalho, and M.~Loss, \textit{ Determination of the spectral
gap for Kac's master equation and related stochastic evolution},
Acta Mathematica {\bf 191}, 1--54 (2003). 



\bibitem{CCL19} E.~ A.~Carlen, M.~C.~Carvalho and M.~P.~Loss, \textit{Chaos, ergodicity and equilibria in a quantum Kac Model}, Adv. in Math., {\bf 358}, 2019, 106827.

\bibitem{CCL3} E.~A.~Carlen, M.~Carvalho and M.~P.~Loss,  \textit{Spectral Gaps for Reversible Markov Processes with Chaotic Invariant Measures: The Kac Process with Hard Sphere Collisions in Three Dimensions},  Ann. of Prob., {\bf 48},  (2020) 2807--2844.


\bibitem{CL93} E.~A.~Carlen and E.~H.~Lieb, \textit{Optimal hypercontractivity for fermi fields and related non-commutative integration inequalities}, Comm. Math. Phys., {\bf 155} (2993) 27--46.


\bibitem{DF90}P.~Diaconis and  J.~Fill, \textit{Strong stationary times via a new form of duality}, Ann. Probab.
{\bf 18} (1990) 1483--1522.

\bibitem{DS93} P.~Diaconis and  L.~Saloff-Coste, \textit{Comparison theorems for reversible Markov chains}, Ann.
Appl. Probab. {\bf 3} (1993), 696--730.


\bibitem{DS81} P.~Diaconis and  M.~Shahshahani, \textit{Generating a random permutation with random transpositions}
Z.Wahrsch. Verw. Gebiete {\bf 57} (1981) 159--179.


\bibitem{Fil16} Y.~Filmus, \textit{Orthogonal basis for functions over a slice of the Boolean hypercube}, Elec. Jour. Comb. {\bf 23} (2016), 1-27. 



\bibitem{FOW19} Y.~Filmus, R.~O'Donnell and X.~Wu, \textit{A Log-Sobolev Inequality for the Multislice, with Applications}, Innovations in Theoretical Computer Science, 2019.



\bibitem{K56} M.~Kac, \textit{Foundations of kinetic theory}, Proc. 3rd Berkeley
symp. Math. Stat. Prob., J. Neyman, ed. Univ. of California, vol 3, (1956)
171--197.

\bibitem{K59}\ M.~Kac \textit{ Probability and Related Topics in Physical Sciences},
Interscience Publ. LTD., London, New York (1959)

\bibitem{K71} K.~Kraus, \textit{General state changes in quantum theory}, Ann. Phys. {\bf 64} (1971), 311--335.

\bibitem{L76} G.~Lindblad, \textit {On the generators of quantum dynamical semigroups}, Comm. Math. Phys., {\bf 48} (1976), 119--130.


%
%
\bibitem{S20} J.~Salez, \textit{A sharp log-Sobolev inequality for the multislice} Annales Henri Lebesgue {\bf 4} (2021) 1143--1161

}
\end{thebibliography}
 \end{document}